# RESEARCH ANNOUNCEMENT

# CRITICAL POINTS ON THE BOUNDARIES OF SIEGEL DISKS

JAMES T. ROGERS, JR.

ABSTRACT. Let $f$ be a polynomial map of the Riemann sphere of degree at least two. We prove that if $f$ has a Siegel disk $G$ on which the rotation number satisfies a diophantine condition, then the boundary of $G$ contains a critical point.

## 1. INTRODUCTION

Let $f : \overline{\mathbf{C}} \to \overline{\mathbf{C}}$ be a polynomial map of the Riemann sphere of degree at least two. Iteration of $f$ leads to a dynamical system and partitions the Riemann sphere into two disjoint sets: the *stable set* or *Fatou set*, and the *unstable set* or *Julia set*. The dynamics of $f$ on the Fatou set is tame, and the dynamics of $f$ on the Julia set is wild.

A point $z$ is a *critical point* if $f'(z) = 0$. A general "fact" about the dynamics of a complex polynomial is that the dynamics of the orbits of the critical points control the dynamics of the map. For instance, the Julia set of $f$ is connected if and only if the orbits of the finite critical points of $f$ are bounded (a subset $A$ of $\overline{\mathbf{C}}$ is *bounded* if the closure of $A$ does not contain infinity). Again the Riemann-Hurwitz relation implies that $f$ has at most $d-1$ finite critical points, where $d$ is the degree of $f$, so that there are bounds on the number of attracting periodic orbits, etc. So it is important to understand the behavior of the orbits of the critical points; among other things, this determines the behavior of the Fatou set.

The work of Sullivan completed the understanding of the dynamics of the polynomial $f$ on the Fatou set. Only a finite number of components of the Fatou set are periodic, and all components of the Fatou set are eventually periodic. The description of the dynamics of $f$ on the Fatou set essentially boils down to describing four possible types of behavior for a component $G$ of the Fatou set satisfying $f(G) = G$.

(1) **Superattracting.** The component $G$ contains a fixed point that is also a critical point.

(2) **Attracting.** The component $G$ contains a fixed point and a critical point, distinct from the fixed point, whose orbit converges to the fixed point.

(3) **Parabolic.** The boundary of $G$ contains a fixed point, and $G$ contains a critical point whose orbit converges to the fixed point.

(4) **Siegel disk.** On the component $G$, $f$ is analytically conjugate to an irrational rotation through the angle $\exp(2\pi i \alpha)$, where $\alpha$ is an irrational real number called the *rotation number* of $G$. The boundary of $G$ is contained in the closure of the orbits of the critical points of $f$.

One natural question concerning components of the Fatou set and critical points of polynomials has remained unanswered:

**Question.** When does the boundary of a Siegel disk contain a critical point?

This question was raised as early as 1982 by Douady [D1, Question 1, p. 40]. See also the survey of Lyubich [L, p. 77], and the work of Ghys [G] and Herman [H1, H2].

It seems reasonable to suspect that the answer is always, but in 1986 M. Herman, building on work of Ghys, proved the startling result that there exist quadratic polynomials with Siegel disks whose boundaries do not contain the critical point. This seems pathological, since it forces the Julia set of that quadratic polynomial to fail to be locally connected.

To avoid such pathology, attention turned to Siegel disks whose rotation numbers $\alpha$ satisfy a diophantine condition: we say $\alpha$ satisfies a diophantine condition if there are numbers $r > 0$ and $k \geq 2$ such that

$$|\alpha - p/q| > r/q^k$$

for every rational number $p/q$. This is a requirement that $\alpha$ be poorly approximated by rational numbers.

Under the assumption that $\alpha$ satisfies a diophantine condition, Ghys [G] showed that the answer is always, provided the boundary of $G$ is a Jordan curve. Let $\varphi : \mathbf{D} \to G$ be a conformal homeomorphism from the unit disk $\mathbf{D}$ to the Siegel disk $G$ that provides the conjugation of $f$ on $G$ to an irrational rotation. If the boundary $B$ of $G$ is a Jordan curve, then $\varphi : \mathbf{D} \to G$ extends to a homeomorphism of $S^1$ onto $B$, and so $f$ is injective when restricted to the boundary $B$. Thus Herman [H1] generalized the Ghys result when he proved the following:

**Theorem** (Herman, 1985). *If the rotation number of $G$ satisfies a diophantine condition, then the boundary of $G$ contains a critical point, provided $f$ is injective when restricted to the boundary of $G$.*

We announce the completion of this program for Siegel disks of complex polynomials by proving the following theorem.

**Theorem 1.** *If the boundary of a Siegel disk $G$ of a polynomial $f$ does not contain a critical point, then $f$ is injective when restricted to the boundary of $G$.*

Thus the following theorem has been proven.

**Theorem 2.** *If the polynomial $f$ has a Siegel disk $G$ and if the rotation number of $G$ satisfies a diophantine condition, then the boundary of $G$ contains a critical point.*

**Comments.** Herman [H1] has proved Theorem 1 for polynomials of the form $z \to z^n + c$. Douady [D2] has an interesting proof of Theorem 1 for quadratic polynomials. The diophantine condition in Theorem 2 is only needed to apply the following theorem of Herman-Yoccoz:

**Theorem.** *If $g$ is an $\mathbf{R}$-analytic diffeomorphism of $S^1$ of rotation number $\alpha$ satisfying a diophantine condition, then $g$ is $\mathbf{R}$-analytically conjugate on $S^1$ to the rotation $z \to \exp(2\pi i\alpha)z$.*

Hence Theorem 2 holds under much weaker forms of the diophantine condition (see [Y]).

## 2. Boundaries of Siegel disks of polynomials

Let $\eta \in S^1 = \partial \mathbf{D}$. The *impression* $I(\eta)$ *of the prime end associated to* $\eta$ is defined by

$$I(\eta) = \{w \in \overline{\mathbf{C}} : \text{there is some half-ray } S = [0, \eta) \text{ in } \mathbf{D}, \text{ perhaps}$$
$$\text{not radial, from 0 to } \eta \text{ such that } w \in \overline{\varphi(S)} - \varphi(S)\}.$$

The proof of Theorem 1 depends on understanding the possibilities for the boundary of a Siegel disk of a polynomial. According to the Structure Theorem of the author [R1, R2], there are two mutually exclusive possibilities for the boundary of such a Siegel disk.

**Structure Theorem.** *The boundary of a Siegel disk $G$ of a polynomial satisfies exactly one of the following*:
  (1) *the impressions of the prime ends of $G$ are disjoint subsets of the boundary of $G$, or*
  (2) *the boundary of $G$ is an indecomposable continuum.*

An *indecomposable continuum* is a compact connected space $X$ that cannot be written as a union $A \cup E$ with $A$ and $E$ connected closed proper subsets of $X$. Indecomposable continua are complicated spaces requiring a special tool kit, so we will deal only with case (1) in this exposition and cover case (2) in [R3].

The Structure Theorem does not give us much information about the topology of the impressions of case (1). The only thing we can determine, aside from the fact that the impressions are disjoint and one-dimensional, is that each impression is *full*, that is, no impression has a bounded complementary domain. This follows from Sullivan's No Wandering Domains Theorem, since a bounded complementary domain of an impression would have to wander under the "irrational rotation" of $f$ (see [R2], p. 817).

Our strategy for exposition will be to present an extremely complicated proof of an extremely simple case of Theorem 1. This proof will convince the reader that other cases follow readily. Hence let $B$ denote the boundary of $G$, and let $x$ and $y$ be two points of $B$ with the property that $f(x) = f(y)$. Let $z$ denote this common value. The point $z$ belongs to exactly one impression $I(\beta)$, so both $x$ and $y$ must belong to the impression $I(\beta - \alpha)$, since $\alpha$ is the rotation number of $G$ and $f^{-1}(I(\beta)) = I(\beta - \alpha)$.

Our extremely simple case is that the impression $I(\beta - \alpha)$ is an arc, a space homeomorphic to $[0, 1]$. Right away the reader can see that since $x$ and $y$ are points on an arc that contain no critical point, the image $I(\beta)$ of this arc must contain a Jordan curve, which contradicts the No Wandering Domains Theorem. But this is not a complicated proof; rather, a complicated proof looks like this:

*Complicated Proof.* Let $K$ be an arc contained in the boundary $B$ of $G$ with the property that $K$ contains two points $x$ and $y$ such that $f(x) = f(y) = z$ and no proper subarc of $K$ contains such a pair of points. In particular, $f$ is injective when restricted to any proper subarc of $K$ and so no proper subarc of $K$ contains both $x$ and $y$. The existence of $K$ follows from an argument using Zorn's Lemma.

Let $Q = \{k \text{ in } K : \text{some point } k' \text{ of } K \text{ distinct from } k \text{ has the property that} f(k) = f(k')\}$. The set $Q$ contains both $x$ and $y$, so it is nonempty. Since $B$ contains no critical point, $Q$ is closed.

Let us assume that $K = A \cup E$, where $A$ and $E$ are arcs. Thus $f$ is injective on $A$ and on $E$, and so their images, $\widetilde{A}$ and $\widetilde{E}$, are also arcs. It follows that $Q \cap A \cap E = \varnothing$. We conclude that if $\widetilde{K}$ denotes the image of $K$, then $\widetilde{K} = \widetilde{A} \cup \widetilde{E}$ is the union of two arcs, and

$$\widetilde{A} \cap \widetilde{E} = \widetilde{A \cap E} \cup \widetilde{Q} \ .$$

Hence $\widetilde{A} \cap \widetilde{E}$ is not connected, and so $H^0(\widetilde{A} \cap \widetilde{E}) \neq 0$. (Our cohomology is reduced Alexander-Čech cohomology with integral coefficients.) From the Mayer-Vietoris sequence of the triple $(\widetilde{K}; \widetilde{A}, \widetilde{E})$ (remember this proof is supposed to be complicated), we see that

$$H^0(\widetilde{A}) \oplus H^0(\widetilde{E}) \to H^0(\widetilde{A} \cap \widetilde{E}) \to H^1(\widetilde{K}) \to H^1(\widetilde{A}) \oplus H^1(\widetilde{E})$$

is an exact sequence. Since $\widetilde{A}$ and $\widetilde{E}$ are arcs, the groups on the ends are trivial, and the middle arrow is thus an isomorphism. Therefore $H^1(\widetilde{K})$ is nontrivial, which is another way of saying that $\widetilde{K}$ has a bounded complementary domain. This contradicts the fact that each impression of $G$ is full.

Let us call a full, one-dimensional compact connected subset of $\mathbf{C}$ a *tree-like continuum*. So each impression of $G$ is a tree-like continuum. If we replace the word "arc" everywhere it occurs in the above proof with the words "tree-like continuum", perhaps the proof will still go through. If so, we have completed the proof of Theorem 1 for case (1). It turns out that the proof goes through, mutatis mutandis, except for one thing, which we invite the reader to find before we reveal it in the next paragraph.

There is one problem: the assumption that $K = A \cup E$ is the assumption that $K$ is a decomposable continuum. We need a more complicated argument for the indecomposable case. We defer this to [R3] with the concluding remark, however, that we have proved in this exposition the following preliminary version of Theorem 2.

**Weak Theorem 2.** *If the polynomial $f$ has a Siegel disk $G$, if the rotation number of $G$ satisfies a diophantine condition, and if the boundary of $G$ does not contain an indecomposable continuum, then the boundary of $G$ contains a critical point.*

## 3. Some topology of the boundary

As a tool to prove this theorem, we obtain the following information on the topology of the boundary of a Siegel disk.

**Theorem 3.** *If $B$ is the boundary of a Siegel disk of a polynomial, then each proper nonempty compact connected subset of $B$ is a tree-like continuum or a point, and any other bounded complementary domain of $B$ is mapped homeomorphically onto $G$ by some iterate of $f$.*

Consequently, while $B$ separates $\mathbf{C}$, no proper closed subset of $B$ separates $\mathbf{C}$. Moreover, if $B$ has more than two complementary domains, then $B$ is, by definition, a "Lakes of Wada continuum". Note there is no assumption on the critical points.

We also show that if the restriction of $f$ to the boundary $B$ of $G$ is injective, then $B$ does not contain a periodic point. In particular, if the boundary of a Siegel disk of a polynomial does not contain a critical point, then the boundary does not contain a periodic point. This gives a partial answer to a question of J. Milnor.

Department of Mathematics, Tulane University, New Orleans, Louisiana 70118
*E-mail address*: `jim@math.tulane.edu`